# A Mixed Integer Programming Model Formulation for Solving the Lot-Sizing Problem

Maryam Mohammadi[1], Masine Md. Tap[2]

[1] Department of Industrial Engineering, Faculty of Mechanical Engineering, Universiti Teknologi Malaysia
Skudai, Johor 81310, Malaysia

[2] Department of Industrial Engineering, Faculty of Mechanical Engineering, Universiti Teknologi Malaysia
Skudai, Johor 81310, Malaysia

**Abstract**
This paper addresses a mixed integer programming (MIP) formulation for the multi-item uncapacitated lot-sizing problem that is inspired from the trailer manufacturer. The proposed MIP model has been utilized to find out the optimum order quantity, optimum order time, and the minimum total cost of purchasing, ordering, and holding over the predefined planning horizon. This problem is known as NP-hard problem. The model was presented in an optimal software form using LINGO 13.0.

## 1. Introduction

Lot-size is defined as the quantity that must be produced or ordered. A manufacturing firm that wants to compete in the market must make the right decisions in lot-sizing problems that have directly effect on the system performance and productivity. This is a critical decision for any manufacturer. The purpose of the lot-sizing problem is to determine the amount of order quantity in each period to meet customer demand over a finite discrete time horizon with the minimum total ordering, purchasing, and holding costs. Smaller lot-size decreases holding cost but raises ordering cost while a larger lot-size reduces ordering cost but leads to the higher inventory costs [1].

A mathematical model that balances holding costs, ordering costs, and purchasing costs must be used to compute optimal or near optimal lot-sizes to minimize cost. According to lean production concepts, small lot-size is preferable, which avoids inventory accumulation, and inventory management and inventory holding costs. Specifically, it would be the lot-size recommended by a mathematical lot-sizing model, which accounts for tradeoff between the associated costs. It must also be taken into account that safety stock is necessary to eliminate the shortage probability in uncertain demand conditions [2].

## 2. Literature Review

Cox and Blackstone [3] defined order management as the scheduling, managing, monitoring and controlling of the operations relevant to customer orders, purchase orders, and production orders. Material requirement planning (MRP) has been utilized to solve the lot-sizing problems with constant demand over a finite time horizon. The economic order quantity (EOQ) has been developed to find the optimum solution while demand remains constant over the time horizon.

The EOQ model was extended by Dye and Hsieh [4] for the variable demand and purchasing cost. The objective was to obtain the optimum replenishment quantity, time planning and cyclic selling price. In order to find the optimal solution, an effective method was proposed using the swarm optimization algorithm. In more general condition, dynamic programming has been recommended. A dynamic programming technique has been developed for single product, multi-period lot-sizing problem [5]. Wang et al. [6] studied the remanufacturing and outsourcing with the single-item and uncapacitated lot-sizing problem. The objective was to determine the lot sizes for manufacturing, remanufacturing, and outsourcing with the minimum costs of holding, setup, and outsourcing. To obtain the optimal solution in the case of existing large amount of returned product, they suggested a dynamic programming technique.

Li, Chena, and Cai [7] considered deterministic time-varying demand, substitutions and return products in the capacitated dynamic lot-sizing problem with batch manufacturing and remanufacturing. To specify the optimal amount of manufactured and remanufactured products in each period, a dynamic programming approach was developed. Hande [8] generalized the well-known previous uncapacitated term and developed polyhedral analysis for the two-item unlimited lot-sizing problem while one-way substitution has been assumed. Finding the minimum cost of production and substitution plan over planning horizon was set as an objective. Li, Chen and Cai [9] considered returned product remanufacturing and demand substitution in the multi-item uncapacitated problem, while disposal and backlog were not allowed. Dynamic programming approach was proposed to minimize the total cost, including manufacturing, remanufacturing, holding and substitution costs.

Even though the above-mentioned techniques find the optimum amount of lot-sizing, high computational resources are required for solving the lot-sizing problem with them [10-13]. Taşgetiren and Liang [14] developed an algorithm based on particle swarm optimization for lot-sizing problem to find order quantities, which will minimize the total ordering and holding costs of ordering decisions. Neural network technique has been used for the single-level lot-sizing problem [15]. Two formulations have been suggested for a stochastic uncapacitated lot-sizing problem [16].

A mixed integer programming model extended the wagner–whitin model. Maximizing the total profit of production and sales over a finite planning horizon was considered as objective while a single-item has been assumed and backordered was not allowed [17]. Absi and Kedad-Sidhoum [18] presented a mixed integer mathematical formulation for solving the capacitated multi product lot-sizing problem when shortage costs and setup times were assumed. They suggested a branch-and-cut framework and fast combinatorial separation algorithms to solve the proposed problem.

The objective of this paper is to develop a model in order to find the optimum order quantity that minimizes the total cost of purchasing, ordering and holding over the specific planning horizon. LINGO software is applied to formulate the proposed mathematical model. A case study is then used to test the model.

## 3. Model Formulation

A mathematical model for lot-sizing problem has been formulated in order to find optimum order quantity and minimize the total purchasing, ordering, and holding costs. A flowchart that represents the model is discussed in the following section.

### 3.1 Flowchart of Optimization Model

There are several factors that need to be taken into consideration in building the model. The factors are lead time to receive parts, the cost of ordering, and holding the parts, required parts in each period, the amount of safety stock for each part, and parts prices in each period. The objective of the model is to determine the cost-effective order quantity that minimizes the total costs of purchasing, ordering, and holding. Figure 1 shows the flowchart for the model formulation. The next section will present the mathematical model formulation.

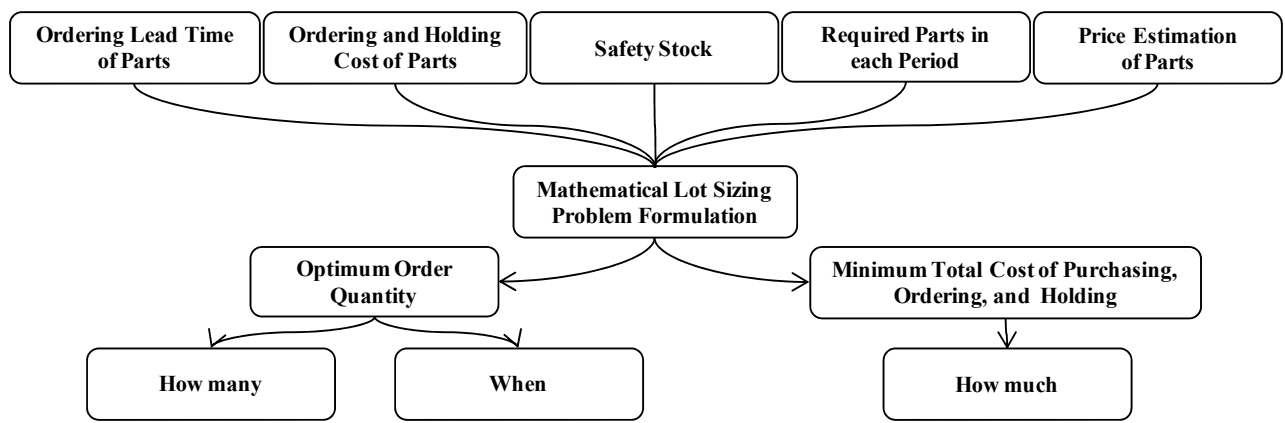

Fig. 1  Model Flowchart

## 3.2 Mathematical Model

There are a few assumptions made in formulating the proposed model. This is because the real situation is dynamic, without the assumption, the model will become too complicated. The assumptions are;
1. Backlogging is not allowed.
2. Number of periods in planning horizon problem is known. Each period is equal to one week and lead time required to receive parts after ordering them is an integer coefficient of a week.
3. Demands for each type of part is estimated and known at the beginning of the planning horizon.
4. Suppliers' lead time is not dependent on volume of ordering size.
5. The price of each part is not constant in entire planning horizon.
6. Total amount of required parts can be ordered over several periods.
7. The ordering cost for parts differ depending on suppliers.

The indices used in the model formulation that indicate the number of particular elements are;
$i$ = index for parts ($i = 1,\ldots,P$)
$t$ = index for periods ($t = 1,\ldots,T$)
$t'$ = index for ordering periods ($t' = 1,\ldots,T$)

List of parameters used in the model formulation are;
$L_i$ = Lead time of part type i
$A_i$ = Ordering cost for part type i
$h_i$ = Holding cost for part type i
$SS_{it}$ = Safety stock required for part type i in period t
$D_{it}$ = Forecasted demand for part type i in period t
$C_{it'}$ = Price of part type i in ordering period t'
$T$ = The number of periods in the planning horizon

Decision variables are as follow;
$Q_{itt'}$ = Amount of required part type i for period t which is ordered in period t'

$Y_{itt'} = \begin{cases} 1 & \text{If part type i required for period t, is ordered in period t',} \\ 0 & \text{Otherwise,} \end{cases}$

Using the above notations, the mathematical model for the lot-sizing problem is presented as follows.

Minimize Total Cost = $f$ (Purchasing Cost + Ordering Cost + Holding Cost)

$$Min\ Z = \sum_{i=1}^{P} \sum_{t=L_i+1}^{T} \sum_{t'=1}^{t-L_i} Q_{itt'} C_{it'}$$
$$+ \sum_{i=1}^{P} \sum_{t=L_i+1}^{T} \sum_{t'=1}^{t-L_i} A_i Y_{itt'}$$
$$+ \sum_{i=1}^{P} \sum_{t=L_i+1}^{T} \sum_{t'=1}^{t-L_i} h_i Q_{itt'}(t-t')$$

(1)

Subject to:

$$\sum_{t'=1}^{t-L_i} Q_{itt'} - SS_{it} \geq D_{it}$$

$i \in \{1,2\ldots P\}, t \in \{L_i+1\ldots T\}$ (2)

$$Q_{itt'} \leq M\ Y_{itt'}$$

$i \in \{1,2\ldots P\}, t \in \{L_i+1\ldots T\}, t' \in \{1\ldots t-L_i\}$ (3)

$$M = Max \sum_{t=1}^{T} D_{it}$$

$i \in \{1,2\ldots P\}$ (4)

$$Y_{itt'} \in \{0,1\}$$

$i \in \{1,2\ldots P\}, t \in \{L_i+1\ldots T\}, t' \in \{1\ldots t-L_i\}$ (5)

$$Q_{itt'} \geq 0$$

$i \in \{1,2\ldots P\}, t \in \{L_i+1\ldots T\}, t' \in \{1\ldots t-L_i\}$ (6)

The objective function presented in Eq. (1) minimizes the total sum of the purchasing cost, the ordering cost, and the holding cost over the planning horizon. The first term shows the purchasing cost of parts over the planning horizon. This cost is obtained by the amount of part type i required for period t ordered in period t' multiply by their associated purchasing costs.

The second term is the ordering cost of parts over the planning horizon. It is the number of times that ordering of part type i have occurred multiply by their associated ordering costs.

The third term calculates the holding cost. It is the amount of part type i required for period t ordered in period t' multiply by the associated cost of inventory holding multiply by (t-t′) that calculates the number of periods in which this ordered amounts were holding.

The first constraint in Eq. (2) ensures that the amount of ordering of each part in previous periods satisfy demand. Furthermore safety stock must be considered in this equation to avoid stock-outs and decrease uncertainty in demand.

The second constraint in Eq. (3) implies that the decision ($Y_{itt'}$) depends on the amount of ($Q_{itt'}$). It means that when there is ordering in period t′, $Y_{itt'}$ is certainly equal to 1 and to satisfy the relation between $Y_{itt'}$ and $Q_{itt'}$, a very large positive number, M is needed.

The third constraint in Eq. (4) suggests one of the feasible values for M. This value must be set sufficiently large to ensure that it is greater than each $Q_{itt'}$. In the absence of backlogging, this amount can be set as the maximum summation of demand for each part over the planning horizon.

The fourth constraint in Eq. (5) indicates that $Y_{itt'}$ is 0 or 1, that is binary decision variable. The last constraint in Eq. (6) is used to define the non-negativity of the decision variable $Q_{itt'}$.

So far, however, there has been little discussion about the multi-item multi-period uncapacitated lot-sizing problem. In addition, the proposed model leads to find the optimum order quantity for each part over several periods by keeping the inventory and ordering costs in a minimum level.

A case study is selected to test the ability of the model, developed to optimize order quantity and minimize total inventory cost.

## 4. Case Study

In order to demonstrate the usability of the suggested model in a dynamic manufacturing environment for multi-period planning horizon, a case study was selected; a trailer manufacturing industry which produces ten different trailers. One of them with twenty-two parts of the trailer has been selected. Table 1 shows the notation given to each part. The planning horizon is established for three months, each month being divided into four periods; therefore, the planning horizon is made of 12 periods.

Table 1: Product Parts

| i = Index for Parts | |
|---|---|
| i = 1 represents part 1 | i = 12 represents part 12 |
| i = 2 represents part 2 | i = 13 represents part 13 |
| i = 3 represents part 3 | i = 14 represents part 14 |
| i = 4 represents part 4 | i = 15 represents part 15 |
| i = 5 represents part 5 | i = 16 represents part 16 |
| i = 6 represents part 6 | i = 17 represents part 17 |
| i = 7 represents part 7 | i = 18 represents part 18 |
| i = 8 represents part 8 | i = 19 represents part 19 |
| i = 9 represents part 9 | i = 20 represents part 20 |
| i = 10 represents part 10 | i = 21 represents part 21 |
| i = 11 represents part 11 | i = 22 represents part 22 |

Referring to Figure 1, one of the first inputs to the model is the ordering lead time of parts. For the case study product, the lead time for ordering each part is shown in Table 2. This is actually the time required from beginning of the purchase request to receive from the outsourcing suppliers.

Table 2: Parts Lead Time (week)

| Part | Lead Time | Part | Lead Time |
|---|---|---|---|
| i = 1 | 1 | i = 12 | 2 |
| i = 2 | 1 | i = 13 | 1 |
| i = 3 | 2 | i = 14 | 1 |
| i = 4 | 1 | i = 15 | 1 |
| i = 5 | 2 | i = 16 | 1 |
| i = 6 | 2 | i = 17 | 1 |
| i = 7 | 2 | i = 18 | 1 |
| i = 8 | 1 | i = 19 | 1 |
| i = 9 | 1 | i = 20 | 1 |
| i = 10 | 1 | i = 21 | 1 |
| i = 11 | 2 | i = 22 | 1 |

The second input is the ordering and holding cost of parts. The ordering and holding costs for the case study product are shown in Tables 3 and 4, respectively.

Table 3: Ordering Cost ($/order)

| Part | Ordering Cost | Part | Ordering Cost |
|---|---|---|---|
| i = 1 | 0.5 | i = 12 | 0.5 |
| i = 2 | 1 | i = 13 | 0.5 |
| i = 3 | 0.5 | i = 14 | 0.5 |
| i = 4 | 0.5 | i = 15 | 0.5 |
| i = 5 | 1 | i = 16 | 1 |
| i = 6 | 1 | i = 17 | 0.5 |
| i = 7 | 0.5 | i = 18 | 1 |
| i = 8 | 1 | i = 19 | 1 |
| i = 9 | 1 | i = 20 | 1 |
| i = 10 | 0.4 | i = 21 | 1 |
| i = 11 | 0.5 | i = 22 | 1 |

Table 4: Weekly Holding Cost ($/unit)

| Part | Holding Cost | Part | Holding Cost |
|---|---|---|---|
| i = 1 | 0.5 | i = 12 | 0.75 |
| i = 2 | 1 | i = 13 | 0.75 |
| i = 3 | 0.5 | i = 14 | 0.5 |
| i = 4 | 0.5 | i = 15 | 0.5 |
| i = 5 | 1 | i = 16 | 1 |
| i = 6 | 1 | i = 17 | 0.75 |
| i = 7 | 0.5 | i = 18 | 1 |
| i = 8 | 1 | i = 19 | 1 |
| i = 9 | 1 | i = 20 | 1 |
| i = 10 | 0.4 | i = 21 | 1 |
| i = 11 | 0.75 | i = 22 | 1 |

The third input to the model is the safety stock. Safety stock for each part is computed using equation (7) that was developed by Walter [19].

Safety Stock = $Z \times \sigma \times \sqrt{LT}$ (7)

It is utilized because previous demand of the trailer is normally distributed. For the case study company, the policy is to set service level at 95 percent, therefore shortage probability is 0.05. According to Normal distribution tables, a probability of 0.05 corresponds to $Z = 1.65$. $\sigma$ is the standard deviation of the past demands of the part. The amount of standard deviation of demand over the last 12 months for each part is shown in Table 5. The lead time (LT) for each part is shown in Table 2.

Table 5: Standard Deviation of Part Demands

| Part | σ | Part | σ |
|---|---|---|---|
| i = 1 | 2.56 | i = 12 | 86.58 |
| i = 2 | 7.28 | i = 13 | 20.96 |
| i = 3 | 7.28 | i = 14 | 2.56 |
| i = 4 | 2.56 | i = 15 | 2.56 |
| i = 5 | 2.56 | i = 16 | 2.56 |
| i = 6 | 2.56 | i = 17 | 2.56 |
| i = 7 | 2.56 | i = 18 | 4.81 |
| i = 8 | 9.67 | i = 19 | 4.81 |
| i = 9 | 2.56 | i = 20 | 14.48 |
| i = 10 | 2.56 | i = 21 | 14.48 |
| i = 11 | 54.42 | i = 22 | 2.56 |

The results of safety stock calculations are presented in Tables 6 and 7. These parts must be held for each period to protect the company against demand fluctuations.

Table 6: Safety Stock (unit)

| Part | Safety Stock | Part | Safety Stock |
|---|---|---|---|
| i = 1 | 4.0 | i = 14 | 4.0 |
| i = 2 | 12.0 | i = 15 | 4.0 |
| i = 3 | 17.0 | i = 16 | 4.0 |
| i = 4 | 4.0 | i = 17 | 4.0 |
| i = 5 | 6.0 | i = 18 | 8.0 |
| i = 6 | 6.0 | i = 19 | 8.0 |
| i = 7 | 6.0 | i = 20 | 24.0 |
| i = 8 | 17.0 | i = 21 | 24.0 |
| i = 9 | 4.0 | i = 22 | 4.0 |
| i = 10 | 4.0 | | |

Table 6 shows the safety stock in terms of unit in each period for parts 1 to 10 and parts 14 to 22.

Table 7: Safety Stock (Kg)

| Part | Safety Stock |
|---|---|
| i = 11 | 126.58 |
| i = 12 | 201.41 |
| i = 13 | 34.03 |

Table 7 shows the safety stock in terms of Kg in each period for parts 11 to 13.

For the fourth input, material requirements planning (MRP) has been used to calculate demand of each part over 12 periods. Results are presented in Tables 8 and 9.

Table 8: Parts Demand (unit/period)

| Period (t) | Part (i) | | | | | | | | | | | | | | | | | | | |
|---|---|---|---|---|---|---|---|---|---|---|---|---|---|---|---|---|---|---|---|---|
| | 1 | 2 | 3 | 4 | 5 | 6 | 7 | 8 | 9 | 10 | 14 | 15 | 16 | 17 | 18 | 19 | 20 | 21 | 22 |
| 1 | -4 | -12 | -17 | -4 | -6 | -6 | -6 | -17 | -4 | -4 | -4 | -4 | -4 | -4 | -8 | -8 | -24 | -24 | -4 |
| 2 | 9 | 35 | -17 | 5 | -6 | -6 | -6 | 13 | 3 | -4 | 3 | -1 | -4 | -4 | -2 | 10 | -24 | -6 | -4 |
| 3 | 14 | 60 | 50 | 20 | 14 | 14 | 13 | 56 | 14 | 11 | 14 | 20 | 5 | 5 | 28 | 30 | 84 | 5 |
| 4 | 16 | 60 | 72 | 20 | 16 | 16 | 20 | 64 | 16 | 16 | 16 | 20 | 16 | 16 | 32 | 32 | 96 | 96 | 16 |
| 5 | 16 | 60 | 72 | 20 | 16 | 16 | 20 | 64 | 16 | 16 | 16 | 20 | 16 | 16 | 32 | 32 | 96 | 96 | 16 |
| 6 | 16 | 60 | 66 | 20 | 16 | 16 | 20 | 64 | 16 | 16 | 16 | 20 | 16 | 16 | 32 | 32 | 96 | 96 | 16 |
| 7 | 16 | 54 | 66 | 18 | 16 | 16 | 18 | 64 | 16 | 16 | 16 | 18 | 16 | 16 | 32 | 32 | 96 | 96 | 16 |
| 8 | 14 | 54 | 66 | 18 | 14 | 14 | 18 | 56 | 14 | 14 | 14 | 18 | 14 | 14 | 28 | 28 | 84 | 84 | 14 |
| 9 | 14 | 54 | 66 | 18 | 14 | 14 | 18 | 56 | 14 | 14 | 14 | 18 | 14 | 14 | 28 | 28 | 84 | 84 | 14 |
| 10 | 14 | 54 | 12 | 18 | 14 | 14 | 18 | 56 | 14 | 14 | 14 | 18 | 14 | 14 | 28 | 28 | 84 | 84 | 14 |
| 11 | 14 | 0 | 0 | 0 | 14 | 14 | 0 | 56 | 14 | 14 | 14 | 0 | 14 | 14 | 28 | 28 | 84 | 84 | 14 |
| 12 | -4 | -12 | -17 | -4 | -6 | -6 | -6 | -17 | -4 | -4 | -4 | -4 | -4 | -4 | -8 | -8 | -24 | -24 | -4 |

Table 8 shows the demand in terms of unit per period for parts 1 to 10 and parts 14 to 22. Negative values indicate that there is no need to order the parts due to adequate availability of stock. As shown in Eq. (2), safety stock must be added to the amount of required part in each period to obtain the final demand per period.

Table 9: Parts Demand (Kg/period)

| Period (t) | Part (i) | | |
|---|---|---|---|
| | 11 | 12 | 13 |
| 1 | -126.58 | -201.41 | -34.03 |
| 2 | -126.58 | -201.41 | -34.03 |
| 3 | 200.66 | 324.22 | 106.66 |
| 4 | 360 | 572.8 | 172 |
| 5 | 360 | 572.8 | 172 |
| 6 | 360 | 572.8 | 172 |
| 7 | 360 | 572.8 | 154.8 |
| 8 | 315 | 501.2 | 154.8 |
| 9 | 315 | 501.2 | 154.8 |
| 10 | 315 | 501.2 | 154.8 |
| 11 | 315 | 501.2 | 0 |
| 12 | -126.58 | -201.41 | -34.03 |

Table 9 shows the demand in terms of Kg per period for parts 11 to 13.

The fifth input to the model is the purchasing price of parts. The purchasing prices of each part over 12 periods are shown in Tables 10 and 11. These are obtained by regression forecasting analyses conducted on historical price of parts.

Table 10: Parts Prices ($/unit)

| Period (t) | Part (i) | | | | | | |
|---|---|---|---|---|---|---|---|
| | 1 | 2 | 3 | 4 | 5 | 6 | 7 |
| 1 | | | | | | | |
| 2 | 380.85 | 4772.27 | 56.68 | 1148.44 | 3324.55 | 182.79 | 91.64 |
| 3 | | | | | | | |
| 4 | | | | | | | |
| 5 | | | | | | | |
| 6 | 381.53 | 4773.39 | 58.35 | 1149.25 | 3326.59 | 183.47 | 92.62 |
| 7 | | | | | | | |
| 8 | | | | | | | |
| 9 | | | | | | | |
| 10 | 382.21 | 4774.51 | 60.02 | 1150.06 | 3328.64 | 184.16 | 93.60 |
| 11 | | | | | | | |
| 12 | | | | | | | |

| Period (t) | Part (i) | | | | | |
|---|---|---|---|---|---|---|
| | 8 | 9 | 10 | 14 | 15 | 16 |
| 1 | | | | | | |
| 2 | 37.06 | 66.30 | 12.51 | 261.27 | 56.68 | 136.18 |
| 3 | | | | | | |
| 4 | | | | | | |
| 5 | | | | | | |
| 6 | 38.09 | 67.63 | 13.40 | 262.28 | 58.35 | 137.09 |
| 7 | | | | | | |
| 8 | | | | | | |
| 9 | | | | | | |
| 10 | 39.12 | 68.96 | 14.29 | 263.28 | 60.02 | 138.01 |
| 11 | | | | | | |
| 12 | | | | | | |

| Period (t) | Part (i) | | | | | |
|---|---|---|---|---|---|---|
| | 17 | 18 | 19 | 20 | 21 | 22 |
| 1 | | | | | | |
| 2 | 172.37 | 16.62 | 24.10 | 130.26 | 143.20 | 91.71 |
| 3 | | | | | | |
| 4 | | | | | | |
| 5 | | | | | | |
| 6 | 175.30 | 17.25 | 24.81 | 133.00 | 145.90 | 92.68 |
| 7 | | | | | | |
| 8 | | | | | | |
| 9 | | | | | | |
| 10 | 178.23 | 17.87 | 25.52 | 135.74 | 148.60 | 93.65 |
| 11 | | | | | | |
| 12 | | | | | | |

Table 11: Parts Prices ($/Kg)

| Period (t) | Part (i) | | |
|---|---|---|---|
| | 11 | 12 | 13 |
| 1 | 1.57 | 0.95 | 0.46 |
| 2 | | | |
| 3 | | | |
| 4 | | | |
| 5 | 1.61 | 0.97 | 0.46 |
| 6 | | | |
| 7 | | | |
| 8 | | | |
| 9 | 1.66 | 0.98 | 0.47 |
| 10 | | | |
| 11 | | | |
| 12 | | | |

## 5. Results and Discussion

The proposed model has been solved using LINGO version 13.0 in personal laptop with 2.76 GHz CPU and 4 GB RAM. Table 12 shows the optimum answer that includes the amount of ordering in 12 periods for part 3, as an example, that have led to the minimum purchasing, ordering, and holding costs while satisfy specified demand completely. In fact, the model gives a plan to the company in order to show how order the parts by satisfying demand and keep the total cost in minimum level.

Required amount of each part per period is indicated as demand ($D_{i,t}$) that has been obtained by MRP calculations. The amount of ordering and ordering period are gained from the model solutions by LINGO 13.0 optimal software to show how many of each part must be ordered per period. This amount is obtained by summation of safety stock and required amount of each part in each period.

The most important finding to appear from the data is that when the cost of purchasing is high, or the ordering is a large amount, the part may be ordered partially in periods that the purchasing cost is lower, although the holding cost has been considered, sometimes this cost has the lower impact on total cost in comparison to purchasing cost impact. Besides that all required amounts of part cannot be ordered in initial periods that have lower price of part, because this affects on the holding cost and total cost will increased considerably. Since LINGO is commanded to minimize the total cost, it makes decision regarding this order. Parts 2, 3, 8, 11, 12, 20, and 21 have been ordered in several periods.

The order quantity (Q) in each period reveals the total amount of each part that must be ordered per period. Q for each period is obtained by summation of the amount of orderings in each order time. For instance, required amount of part 3 in period 3 is ordered in period 1, and some of the demand of it in periods 4, 5, 6 and 7 are also ordered in period 1; therefore, the order quantity for period 1 will be sum of all ordering amounts in period 1. So, the company must order the total amount of 171 of part 3 in period 1. Accordingly, the order quantity for the other parts and periods are obtained by this manner. Number of ordering shows how many times the ordering per period has been placed. Similar calculations have been done for the other parts.

Table 12: Model Optimum Solution

| Part (i) | Period (t) | Demand ($D_{i,t}$) | Amount of Ordering (Unit) | Ordering Period (t′) | Order Quantity (Q) | Number of Ordering | Total Cost (Z) |
|---|---|---|---|---|---|---|---|
| 3 | 1 | -17 | 0 | - | 171 | 5 | 36523.68 |
| | 2 | -17 | 0 | - | 100 | 3 | |
| | 3 | 50 | 67 | 1 | 160 | 4 | |
| | 4 | 72 | 89 | 1 | 80 | 2 | |
| | 5 | 72 | 9,40,40 | 1,2,3 | 20 | 1 | |
| | 6 | 66 | 3,40,40 | 1,2,3 | 23 | 1 | |
| | 7 | 66 | 3,20,20,20,20 | 1,2,3,4,5 | 23 | 1 | |
| | 8 | 66 | 60,23 | 3,6 | 46 | 2 | |
| | 9 | 66 | 60,23 | 4,7 | 0 | 0 | |
| | 10 | 12 | 29 | 8 | 0 | 0 | |
| | 11 | 0 | 17 | 8 | 0 | 0 | |
| | 12 | -17 | 0 | 0 | 0 | 0 | |

The minimum total cost incurred to purchase, order, and hold each part over 12 periods is shown in Table 13. The total cost for all parts is $ 4,423,918.29.

Table 13: Total Cost of Optimum Solution for each Part over Planning Horizon ($)

| Part | Total Cost (Z) | Part | Total Cost (Z) |
|---|---|---|---|
| i = 1 | 69885.85 | i = 12 | 17270.80 |
| i = 2 | 2917218.12 | i = 13 | 1879.97 |
| i = 3 | 36523.68 | i = 14 | 46468.01 |
| i = 4 | 226446.94 | i = 15 | 11088.54 |
| i = 5 | 625584.80 | i = 16 | 22230.02 |
| i = 6 | 34854.82 | i = 17 | 28161.04 |
| i = 7 | 18531.30 | i = 18 | 6298.28 |
| i = 8 | 28388.13 | i = 19 | 9223.96 |
| i = 9 | 12106.40 | i = 20 | 128969.36 |
| i = 10 | 2269.25 | i = 21 | 151998.00 |
| i = 11 | 13442.16 | i = 22 | 15078.87 |

One of the significant findings which emerges from this study is to attain optimum order quantity in each period and minimum total cost of purchasing, ordering, and holding the parts.

## 6. Conclusions

In this paper, a mixed integer programming (MIP) model formulation has been created. The case study company needs safety stock in order to prevent stock-out, hence one method based on the variance in demand, lead time and the target service level is used to determine the amount of safety stock for each part. MRP is done to obtain required amount of each part per period. The parts prices are estimated for 12 periods using regression forecasting method.

The assumptions, parameters, and the decision variables have been defined. Then, the mixed integer model formulation included the objective function, and constraints have been demonstrated and discussed. The model is solved using LINGO 13.0. Results show that the model constraints are satisfied as all demands of parts are planned to be ordered in a correct amount when they are needed. The proposed MIP model is able to attain optimum order quantity in each period and minimum total cost of purchasing, ordering, and holding the parts.

Metaheuristic algorithms such as genetic algorithm and tabu search can also be used to solve the proposed MIP model. Comparison can then be made to study the best method.

**Acknowledgments**

The author would like to thank the Case Study Company for their assistance and support of this work. The author would also like to thank Universiti Teknology Malaysia for their support that has made this research work possible.

**Maryam Mohammadi,** received Bachelor's Degree in Industrial Engineering (Systems Planning and Analysis) in Islamic Azad University - South Tehran Branch in 2009, and Master's degree in Industrial Engineering program from Universiti Teknologi Malaysia (UTM) in 2012. She is currently a fulltime research student of Universiti Teknologi Malaysia (UTM) enrolled in the Doctor of Philosophy (Mechanical Engineering) program. Her research interests are in production planning and control, supply chain management, operations research, optimization, and simulation.

**Masine Md. Tap,** received her Bachelor's Degree in Mechanical Engineering from Universiti Teknologi Malaysia (UTM) in 1986, MPhil in Computer Aided Engineering from Herriot-Watt University, United Kingdom in 1989 and PhD. from Dundee University, United Kingdom in 1999. She is now an associate professor and the course coordinator for MEng. (Industrial Engineering) program in faculty of Mechanical Engineering in Universiti Teknologi Malaysia (UTM).